

\input amstex

\define\scrO{\Cal O}
\define\Pee{{\Bbb P}}
\define\Zee{{\Bbb Z}}

\define\Ar{{\Bbb R}}

\define\proof{\demo{Proof}}
\define\endproof{\qed\enddemo}
\define\endstatement{\endproclaim}
\define\theorem#1{\proclaim{Theorem #1}}
\define\lemma#1{\proclaim{Lemma #1}}
\define\proposition#1{\proclaim{Proposition #1}}
\define\corollary#1{\proclaim{Corollary #1}}
\define\claim#1{\proclaim{Claim #1}}

\define\section#1{\specialhead #1 \endspecialhead}
\define\ssection#1{\medskip\noindent{\bf #1}}

\documentstyle{amsppt}
\pageno=1
\topmatter
\title Vector bundles  and $SO(3)$-invariants for elliptic surfaces I
\endtitle
\author Robert Friedman
\endauthor
\address Columbia University, New York, NY 10027
\endaddress
\email rf\@math.columbia.edu \endemail
\thanks Research partially supported by NSF grant  DMS-9203940
\endthanks
\subjclass Primary  14J60, 57R55; Secondary 14D20, 14F05, 14J27
\endsubjclass
\endtopmatter
\input cyracc.def
\newfam\cyrfam
\font\tencyr=wncyr10
\font\sevencyr=wncyr7
\font\fivecyr=wncyr5
\def\cyr{\fam\cyrfam\tencyr\cyracc}
\textfont\cyrfam=\tencyr
\scriptfont\cyrfam=\sevencyr
\scriptscriptfont\cyrfam=\fivecyr

\def\Shah{{\cyr Sh}}
\document

\section{1. Introduction.}

Beginning with Donaldson's seminal paper on the failure of the $h$-cobordism
theorem in dimension 4 [4], the techniques of gauge theory have proved to be
highly successful in analyzing the smooth structure of simply connected
elliptic
surfaces. Recall that a simply connected elliptic surface $S$ is specified up
to
deformation type by its geometric genus $p_g(S)$ and by two relatively prime
integers $1\leq m_1\leq m_2$, the multiplicities of its multiple fibers.
Here, if $p_g(S) = 0$, a surface $S$ such that $m_1 =1$ is rational, and thus
all surfaces $S$ with $p_g(S) = 0$ and $m_1=1$ are deformation equivalent and
in particular diffeomorphic. Moreover, if $p_g(S) = 1$ and $m_1 = m_2 =1$,
then $S$ is a $K3$ surface. In all other cases, $S$ is a surface with Kodaira
dimension one.

Our goal in this series of three papers is to prove the following result:

\theorem{} Two possibly blown up simply connected elliptic surfaces are
diffeomorphic if and only if they are deformation equivalent. More precisely,
suppose that $S$ and $S'$ are relatively minimal simply connected elliptic
surfaces. Suppose that $S$ has multiple fibers of multiplicities $m_1$ and
$m_2$, with $1\leq m_1 \leq m_2$, and that $S'$ has multiple fibers of
multiplicities $m_1'$ and $m_2'$, with $1\leq m_1' \leq m_2'$. Let $\tilde S$
be a blowup of $S$ at $r$ points and $\tilde S'$ a blowup of $S'$ at $r'$
points. Suppose that $\tilde S$ and $\tilde S '$ are diffeomorphic. Then $r=r'$
and $p_g(S) = p_g(S')$, and moreover: \roster
\item"{(i)}" If $p_g(S)> 0$, then $m_1=m_1'$ and $m_2=m_2'$.
\item"{(ii)}" If $p_g(S) = 0$, then $S$ is rational, i.e\. $m_1=1$, if and only
if $S'$ is rational if and only if $m_1'=1$. If $S$ and $S'$ are not rational,
then $m_1=m_1'$ and $m_2=m_2'$.
\endroster
\endstatement
\smallskip

There is also a routine generalization to the case of a finite cyclic
fundamental group. The statements in the theorem that $r=r'$ and $p_g(S) =
p_g(S')$ are easy consequences of the fact that $\tilde S$ and $\tilde S'$ are
homotopy equivalent, and the main point is to determine the multiplicities.
Before discussing the proof of the theorem in more detail, we shall review some
of the history of the classification of simply connected elliptic surfaces:

\theorem{1.1 [8]} There is a function $f(m_1, m_2)$ defined on pairs of
relatively prime positive integers $(m_1, m_2)$ such that $f$ is symmetric and
finite-to-one provided that neither $m_1$ nor $m_2$ is $1$, with the following
property: Let $S$ and $S'$ be two simply connected surfaces with $p_g(S) = 0$.
Denote the multiplicities of the multiple fibers of $S$ by $m_1, m_2$ and the
multiplicities for $S'$ by $m_1', m_2'$. If $S$ and $S'$ are diffeomorphic,
then
$f(m_1, m_2) = f(m_1', m_2')$. Moreover, let $\tilde S$ and $\tilde S'$ be
blowups of $S$ and $S'$ at $r$ points. Then:
\roster
\item"{(i)}" Every diffeomorphism $\psi \: \tilde S \to \tilde S'$ pulls back
the cohomology class of an exceptional curve on $\tilde S'$ to $\pm$ the
cohomology class of an exceptional curve on $\tilde S$.
\item"{(ii)}" Every diffeomorphism $\psi \: \tilde S \to \tilde S'$ pulls back
the cohomology class of a general fiber on $\tilde S'$ to a rational multiple
of the
cohomology class of a  general fiber on $\tilde S$.
\item"{(iii)}" If $\tilde S$ and $\tilde S'$ are diffeomorphic, then  $f(m_1,
m_2) = f(m_1', m_2')$. \qed
\endroster
\endstatement

The function $f(m_1, m_2)$ was then determined by S\. Bauer [1] (the case
$m_1 = 2$ is also in [8]):

\theorem{1.2} In the above notation,
$$f(m_1, m_2) = \frac{(m_1^2-1)(m_2^2-1)}3 -1.\qed$$
\endstatement

For the case $p_g(S) > 0$, there is the following result [9]:

\theorem{1.3}  Let $S$ and $S'$ be two simply connected surfaces with
$p_g(S) > 0$. Denote the multiplicities of the multiple fibers of $S$ by $m_1,
m_2$ and the multiplicities for $S'$ by $m_1', m_2'$. If $S$ and $S'$ are
diffeomorphic, then $m_1m_2 = m_1'm_2'$. Moreover, let $\tilde S$ and $\tilde
S'$ be blowups of $S$ and $S'$ at $r$ points. Then:
\roster
\item"{(i)}" Every diffeomorphism $\psi \: \tilde S \to \tilde S'$ pulls back
the cohomology class of an exceptional curve on $\tilde S'$ to $\pm$ the
cohomology class of an exceptional curve on $\tilde S$.
\item"{(ii)}" Except possibly for $p_g(S) =1$, every diffeomorphism $\psi \:
\tilde S \to \tilde S'$ pulls back the cohomology class of a general fiber on
$\tilde S'$ to a rational multiple of the
cohomology class of a  general fiber on $\tilde S$.
\item"{(iii)}" If $\tilde S$ and $\tilde S'$ are diffeomorphic, then
$m_1m_2 = m_1'm_2'$. \qed
\endroster
\endstatement

The crux of the argument involves calculating a coefficient of a suitable
Donaldson polynomial invariant $\gamma _c(S)$. In fact, it is shown in [9] that
$\gamma _c(S)$ can be written as a polynomial in the intersection form $q_S$
and the primitive class $\kappa$ such that the class of a general fiber  $[f]$
of $S$ is equal to $m_1m_2\kappa$, and that, for $c$ sufficiently large, the
first nonzero coefficient of this polynomial is given as follows: let $n =
2c-2p_g(S) -1$ and $d = 4c - 3p_g(S) -3$. If $\gamma _c(S) = \sum
_{i=0}^{[d/2]}a_iq_S^i\kappa ^{d-2i}$, then $a_i = 0$ for $i>n$ and
$$a_n = \frac{d!}{2^nn!}(m_1m_2)^{p_g(S)}.$$
The proof of this statement involves showing that the moduli space of stable
vector bundles $V$ with $c_1(V) = 0$ and $c_2(V) = c$ fibers holomorphically
over a Zariski open subset of a projective space, and that the fiber consists
of $m_1m_2$ copies of a complex torus.  It is natural to wonder if the
techniques
of [9] can be pushed to determine some of the remaining terms. However, it
seems
to be difficult to use the vector bundle methods used in [9] to make the
necessary calculations, even in the case of no multiple fibers. Thus, it is
natural to look for other techniques to complete the $C^\infty$ classification
of elliptic surfaces.

Using a detailed analysis of certain moduli spaces of vector bundles, Morgan
and O'Grady [15] together with Bauer [2] were  able to calculate the
coefficient $a_{n-1}$ in case $p_g(S) = 1$ and $c=3$. The calculation is long
and involved for the following reason: the moduli spaces are nonreduced, not
necessarily of the correct dimension, and (in the case of trivial determinant)
the integer $c$ is not in the ``stable range." The final answer is that, up to
a
universal combinatorial factor, $a_{n-1} = m_1m_2(2m_1^2m_2^2 - m_1^2 -
m_2^2)$.
{}From this and from the knowledge of $m_1m_2$, it is easy to determine the
unordered pair $\{m_1, m_2\}$. In addition, the calculation shows that the
class
of a fiber of $S$ is preserved up to rational multiples in case $p_g(S) = 1$ as
well (the possible exception in (ii) of Theorem 1.3 above), provided that not
both of $m_1$ and $m_2$ are 1.

In  the
proof of the main theorem, we shall use the following results.  Aside from
standard techniques in the theory of vector bundles, and the gauge theory
results that are described in the book [9], we use only the results of this
series of papers and of [9] to handle the case $p_g(S) >0$. In case $p_g(S) =
0$, we use the results in this series and in [8], as well as the calculation of
Bauer described  in Theorem 1.2 in case $m_1m_2 \equiv
0\mod 2$. In case $m_1m_2\equiv 1 \mod 2$, our proof does not depend on Bauer's
results.

Next we outline the strategy of the argument. Following a well-established
principle [6], [10], we shall work with $SO(3)$-invariants instead
of $SU(2)$-invariants since these are often much easier to calculate.
Moreover, in case $b_2^+=1$ a good choice of an $SO(3)$ invariant can
simplify the problem that the invariant depends on the choice of a certain
chamber. Thus we must choose a class $w\in H^2(S; \Zee/2\Zee)$ to be the second
Stiefel-Whitney class of a principal $SO(3)$ bundle, although it will usually
be more convenient ot work with a lift of $w$ to $\Delta \in H^2(S; \Zee)$. One
possible choice of a lift of $w$ would be the class $\kappa$, the primitive
generator of $\Zee ^+\cdot [f]$, or perhaps $c_1(S)$, or even $[f]$. All of
these classes are rational multiples of $[f]$, and they do not simplify the
problem.

Instead we shall consider the case where $\Delta$ is transverse to $f$, more
specifically where $\Delta \cdot \kappa = 1$. Of course, we shall need to
choose $\Delta$ to be the class of a holomorphic divisor as well in order to be
able to apply algebraic geometry. As we shall see in Section 2, we can always
make the necessary choices and the final calculation will show that the answer
does not depend on the choices made. Note that $\Delta$ is well-defined up to a
multiple of $\kappa$, and that the choices $\Delta$ and $\Delta -\kappa$
correspond to different choices for $w^2\equiv p\mod 4$. Finally, as we shall
show in Section 3, in case $b_2^+(S) =1$ or equivalently $p_g(S)=0$, there is a
special chamber $\Cal C(w,p)$ which is natural in an appropriate sense under
diffeomorphisms.

With this choice of $\Delta$, the study of the relevant vector bundles divides
into two very different cases, depending on whether $m_1m_2\equiv 0\mod 2$ or
$m_1m_2\equiv 1\mod 2$. In this paper we shall collect results which are needed
for both cases and show how the main theorems follow from the calculations in
Parts II and III. In Part II of this series, we shall consider the case where
$\Delta \cdot \kappa = 1$ and $m_1m_2\equiv 0\mod 2$. In this case, $m_1$, say,
is even. Since $\Delta  \cdot f =m_1m_2 \Delta \cdot \kappa$, a vector bundle
$V$ with $c_1(V) = \Delta$ has even  degree on a general fiber $f$. At first
glance, then, it seems as if we are again in the situation of [7] and [9]
and that there is no new information to be gained from the Donaldson
polynomial. However, it turns out that the asymmetry between $m_1$ and $m_2$
appears in the moduli space as well. In this case, the moduli space again
fibers
holomorphically over a Zariski open subset of a projective space. But the
fibers
now consist of just $m_2$ copies of a complex torus. We then have, by an
analysis that closely parallels [9], the following result:

\theorem{1.4} Let $w$ and $p$ be as above, and set $d = -p - 3(p_g(S) + 1)$ and
$n = (d-p_g(S))/2$. Suppose that $\gamma _{w,p}(S)$ is
the Donaldson polynomial for the $SO(3)$-bundle $P$ over $S$ with $w_2(P) = w$
and $p_1(P) = p$ where if $p_g(S) =0$ this polynomial is associated to the
chamber $\Cal C(w,p)$ defined in \rom{(3.6)} below. Then, assuming that $m_1$
is even and writing  $\gamma _{w,
p}(S)= \sum _{i=0}^{[d/2]}a_iq_S^i\kappa _S^{d-2i}$, we have, for all
$p$ such that $-p\geq 2(4p_g(S) + 2)$, $a_i= 0$ for $i>n$ and  $$a_n =
\frac{d!}{2^nn!}(m_1m_2)^{p_g(S)}m_2.$$
\endstatement
\medskip

In particular, the leading coefficient contains an ``extra" factor of
$m_2$. Using this and either [9] in case $p_g(S) > 0$ or [1] in case $p_g(S)
= 0$, we may then determine $\{m_1, m_2\}$. Note that in case $p_g(S) = 0$ and
one of $m_1$, $m_2$ is 1, then it cannot be $m_1$ since $m_1$ is even. Thus
$m_2
=1$ and the leading coefficient does not determine $m_1$ (as well it cannot).

Finally, in Part III we shall discuss the case where $\Delta \cdot \kappa = 1$
and $m_1m_2\equiv 1\mod 2$.
If $m_1m_2 \equiv 1 \mod 2$, then vector bundles
$V$ with $c_1(V ) = \Delta$ have odd degree when restricted to a general fiber,
and the general methods for studying vector bundles on elliptic surfaces
described in [7] and [9] Chapter 7 do not apply. Thus we must develop new
techniques for studying such bundles, and this is the subject of Part III of
this series. Fortunately, it turns  out that this moduli problem is in many
ways much simpler to study than the case of even degree on the general fiber.
For example, as long as the expected dimension is nonnegative, for a suitable
choice of ample line bundle the moduli space is always nonempty, irreducible,
and smooth of the expected dimension. Moreover a Zariski open subset of the
moduli space is independent of the multiplicities, and from this one can show
easily that the leading coefficient of the Donaldson polynomial for the
corresponding $SO(3)$-bundle is (up to the usual combinatorial factors) equal
to 1. At first glance, this rather disappointing result suggests that no new
information can easily be gleaned from the Donaldson polynomial. However, this
suggestion is misleading: in some sense, the structure of the moduli space
allows the contribution of the multiple fibers to be localized around the
multiple fibers, enabling us to calculate the next two coefficients in the
Donaldson polynomial. By contrast, in the case of trivial determinant,
the moduli space for a surface with two multiple fibers of multiplicities $m_1$
and $m_2$ looks roughly like a branched cover of the corresponding moduli space
for a surface without multiple fibers. A further simplification is that we can
work with moduli spaces of small dimension, for example dimension two or four.
Using the vector bundle results, we shall show:

\theorem{1.5} Let $S$ be a simply connected elliptic surface with two multiple
fibers $m_1$ and $m_2$, with $m_1m_2 \equiv 1\mod 2$. Let $w \in H^2(S;
\Zee/2\Zee)$ satisfy $w\cdot \kappa = 1$. Suppose that $\gamma _{w,p}(S)$ is
the Donaldson polynomial for the $SO(3)$-bundle $P$ over $S$ with $w_2(P) = w$
and $p_1(P) = p$ where if $p_g(S) =0$ this polynomial is associated to the
chamber $\Cal C(w,p)$ defined in \rom{(3.6)} below.
\roster
\item"{(i)}" Suppose $w$ and $p$ are chosen so that the expected complex
dimension of the moduli space $-p -3(p_g(S) +1)$ is $2$. Then for all $\Sigma
\in H_2(S; \Zee)$,
$$\gamma _{w,p}(S)(\Sigma, \Sigma) = \Sigma ^2 + ((m_1^2m_2^2)(p_g(S) + 1)-
m_1^2-m_2^2) (\Sigma \cdot \kappa)^2.$$
\item"{(ii)}" Suppose $w$ and $p$ are chosen so that the expected complex
dimension of the moduli space $-p -3(p_g(S) +1)$ is $4$. Then for all $\Sigma
\in H_2(S; \Zee)$,
$$\gamma _{w,p}(S)(\Sigma, \Sigma ,\Sigma, \Sigma ) = 3(\Sigma ^2)^2 +
6C_1(\Sigma ^2 )(\Sigma \cdot \kappa)^2 + (3C_1^2 -
2C_2)(\Sigma \cdot
\kappa)^4, $$
where
$$\align
C_1&= (m_1^2m_2^2)(p_g(S) + 1)-m_1^2-m_2^2;\\
C_2 &= (m_1^4m_2^4) (p_g(S) + 1) -m_1^4-m_2^4.
\endalign$$
Here $C_1$ is the second coefficient of the degree two polynomial.
\endroster
\endstatement
\medskip

Note that the final answer has the following self-checking features. First, it
is a polynomial in $q_S$ and $\kappa _S$. If $p_g(S) = 1$ and $m_1 = m_2 = 1$,
so
that $S$ is a $K3$ surface, then the term $(\Sigma \cdot \kappa)$ does not
appear. This is in agreement with the general result that $\gamma _{w,p}(S)$ is
a
multiple of a power of $q_S$ alone. If $p_g(S) = 0$ and $m_1 =1$, then the
answer is independent of $m_2$, since in this case all of the surfaces $S$ for
various choices of $m_2$ are diffeomorphic. In fact, we shall turn this remark
around and use the knowledge of  $\gamma _{w,p}(S)$ for $p_g =0$, $m_2
=1$ and $m_1$ arbitrary, to determine $\gamma _{w,p}(S)$ in general.

The techniques used to prove Theorem 1.5 should be capable
of further generalization. For example, these methods should give in principle
(that is, up to the knowledge of the multiplication table for divisors in
$\operatorname{Hilb}^nS$), the full polynomial invariant in case $m_1$ and
$m_2$ are odd. One might make a conjectural formula for  $\gamma
_{w,p}(S)$ in general along the lines suggested by Kronheimer and Mrowka in
[13]. In our case the formula should conjecturally read as follows: let $\gamma
_t(\Sigma)$ be the Donaldson polynomial $\gamma _{w,p}(S)(\Sigma , \dots ,
\Sigma)$ for $w = \Delta \mod 2$ or $w= \Delta -\kappa \mod 2$ and $p$ chosen
so that $w^2 \equiv p \mod 4$ and $-p - 3\chi (\scrO _S) =2t$, so that the
complex dimension of the moduli space is $2t$. It follows from Proposition
2.1 below that $\gamma _t$ depends only on $t$. Then the natural analogue of
the conjectures in [13] is the conjecture that
$$\sum _{t\geq 0}\frac{\gamma _t(\Sigma)}{(2t)!}
= \exp \fracwithdelims(){q_S}2\frac{\bigl(\cosh ( m_1m_2(\kappa \cdot
\Sigma))\bigr)^{p_g+1}}{\cosh (m_1(\kappa \cdot \Sigma)) \cosh (m_2(\kappa
\cdot
\Sigma))}.$$
It essentially follows from Theorem 1.5 that this formula is
correct through the first three terms, including the case $p_g =0$ where the
quotient is not given by a finite sum of exponentials, and it is likely that a
further extension of the methods in Part III and the knowledge of the
multiplication in  $\operatorname{Hilb}^nS$ can establish the general formula.
Finally we should add that many of these techniques  have  applications to the
$SU(2)$ case.

Morgan and Mrowka [14] have independently determined the second coefficient
$a_{n-1}$ for all $S$ such that $p_g(S) \geq 1$, for the case of the
$SU(2)$-invariant $\gamma _c(S)$. The answer is that, up
to combinatorial factors,  $$a_{n-1} = (m_1m_2)^{p_g(S)}((m_1^2m_2^2)(p_g(S) +
1)-m_1^2-m_2^2).$$
{}From this, it is again easy to see that the diffeomorphism type of
$S$ determines the unordered
pair $\{m_1, m_2\}$ in case $p_g(S) \geq 1$. The proof of this formula uses the
knowledge of $a_{n-1}$ for the case of $p_g(S)=1$, together with the gauge
theory gluing techniques developed by Mrowka in [16], to determine the
coefficient  $a_{n-1}$ for $p_g(S) > 1$.

Let us finally describe the contents of this paper. In Section 2 we discuss the
possible choices for $w=w_2(P)$ up to diffeomorphisms of $S$ and show that
there
is always a generic elliptic surface for which $w$ is the mod 2 reduction of a
holomorphic divisor. There is also a discussion of certain elliptic surfaces
which can be constructed from $S$. In Section 3 we introduce a class of ample
line bundles which we shall use to define stability and which are well-adapted
to the geometry of $S$. Section 4 explains the meaning of stability of a vector
bundle $V$ with respect to such a line bundle: stability is equivalent to the
assumption that the restriction of $V$ to almost every fiber is semistable.
Finally, in Section 5, we show how  the main results concerning Donaldson
polynomials lead to $C^\infty$ classification results.

\section{2. Preliminaries on elliptic surfaces.}

Let $S$ be a simply connected elliptic surface with at most two multiple fibers
of multiplicities $m_1 \leq m_2$. Here we shall allow $m_1$ or both $m_1$ and
$m_2$ to be one. Let $[f]$ denote the class in homology of a smooth
nonmultiple fiber of $S$. There is a unique homology class $\kappa _S =
\kappa$ such that $[f] = m_1m_2\kappa$, and $\kappa$ is primitive [9]. Let
$P$  be a principal $SO(3)$-bundle over $S$ with $w_2(P) = w$ and $p_1(P) = p$.
Note that $w^2 \equiv p \mod 4$.
We shall be concerned with bundles $P$ such that $w\cdot \kappa\bmod 2 = 1$. In
this section, we shall show that, modulo diffeomorphism, the choice of $w$ is
not
essential. Indeed, we shall prove that, given a class
$w$, there is a diffeomorphism $\psi \: S \to S'$, where $S'$ is again a simply
connected elliptic surface with two multiple fibers of
multiplicities $m_1$ and $m_2$, such that $\psi ^*\kappa _{S'} = \kappa _S$ and
such that there exists a holomorphic divisor $\Delta$ with $w=\psi ^*[\Delta]
\mod 2$. Thus we may always assume that $w$ is the reduction of a $(1,1)$
class.
We begin with an arithmetic result, which is not in fact needed in what follows
but which helps to clarify the role of the choice of $w$ modulo
diffeomorphisms.
In the arguments below, we shall sometimes blur the distinction between
$H_2(S)$
and $H^2(S)$ using the canonical identification between these two groups.

\proposition{2.1} Let $S$ be a simply connected elliptic surface.
\roster
\item"{(i)}" Suppose that $m_1m_2 \equiv 1\mod 2$, and let $a\in \Zee/4\Zee$.
Then the group of
orientation-preserving diffeomorphisms $\psi \: S \to S$  such that
$\psi _*([f]) = [f]$ acts transitively on the set of $w\in H^2(S; \Zee/2\Zee)$
such that $w\cdot \kappa =1$ and $w^2 \equiv a \mod 4$.
\item"{(ii)}" If $m_1m_2 \equiv 0\mod 2$ and $a\in \Zee/4\Zee$, then there are
at
most three orbits of the set
$$\{\, w\in H^2(S; \Zee /2\Zee): w\cdot \kappa =1 \text{ and } w^2\equiv
a \mod 4\,\}$$
under the group of diffeomorphisms of $S$
which fix $\kappa$.
\endroster
\endstatement
\proof   Let $L$ be
the image of $H_2(S - \pi ^{-1}(D))$ in $H_2(S)$, where $D$ is a small disk
in $\Pee ^1$ which we may assume contains the multiple fibers and no other
singular fiber. Thus $L\subseteq (\kappa^\perp)$, and in fact $L$ has index
$m_1m_2$ in $(\kappa ^\perp)$. Let $\varphi$ be an automorphism of the
lattice $H_2(S; \Zee)$ fixing $\kappa$. Thus by restriction $\varphi$ induces
an automorphism of $(\kappa ^\perp)$. The method of proof of Theorem 6.5 of
Chapter 2 of [9] shows that there is a diffeomorphism $\psi$, automatically
orientation-preserving,  inducing $\varphi$ provided that $\varphi(L) \subseteq
L$ and that $\varphi$ has real spinor norm one.

Clearly we may write $L = \Zee [f] \oplus W$, where $W$ is an even unimodular
lattice. Moreover $(\kappa ^\perp) = \Zee \cdot\kappa \oplus W$, with the
inclusion $L\subseteq (\kappa ^\perp)$ the natural inclusion given by $[f] =
m_1m_2\kappa$. If $W^\perp$ denotes the orthogonal complement of $W$ in
$H_2(S; \Zee)$, then $W^\perp = \operatorname{span}\{\kappa, x\}$ for some
class $x$ with $x\cdot \kappa = 1$. Given $a\bmod 4$, we can always assume
after replacing $x$ by $x+\kappa$ that $x^2 \equiv a \mod 4$. Now it is easy
to describe all automorphisms of $H_2(S; \Zee)$ fixing $\kappa$: choosing an
isometry $\tau$ of $W$, $\varphi $ is given by
$$\align
\varphi (\kappa ) &= \kappa;\\
\varphi (\alpha) &=  \tau(\alpha) +\ell(\alpha)\kappa, \quad\alpha \in W;\\
\varphi (x) &= x + c\kappa + \beta.
\endalign$$
Here $\ell$ is an arbitrary homomorphism $W\to \Zee$ and $\beta$ is the
unique element of the unimodular lattice $W$ such that $-\beta \cdot \alpha =
\ell (\tau (\alpha))$ for all $\alpha \in W$. Furthermore $c=-\beta ^2/2$. It
is clear that every choice of $\tau$ and $\ell$ (or equivalently $\beta$)
produces an automorphism $\varphi$, and that $\varphi(L) = L$ if and only if
$m_1m_2$ divides $\ell$ or equivalently $\beta$.
If  $x'$ is another class such that $x'\cdot \kappa \equiv 1\mod 2$ and
$(x')^2 \equiv x^2 \mod 4$, we can write $x' = nx + b\kappa + \beta$, where
$\beta \in W$. Since we only care about $x'\bmod 2$, we may assume that $n=1$.
Note that $2b + \beta ^2 \equiv 0\mod 4$ and thus $b \equiv \beta ^2/2 \mod 2$.

First assume that $m_1m_2$ is odd. Then since $\beta \equiv m_1m_2\beta \mod
2$, we may assume that $\beta $ is divisible by $m_1m_2$. Choosing $\tau =
\operatorname{Id}$ and $\ell$, $c$ in the definition of $\varphi$ as specified
by $\beta$ gives $\varphi$ such that $\varphi (x) \equiv x' \mod 2$. As
$\varphi$ is unipotent, it is easy to see that $\varphi$ has spinor norm one,
i.e\. that $\varphi$ is in the same connected component of the group of
automorphisms of the quadratic form of $H_2(S; \Ar)$ as the identity. Thus
there is a diffeomorphism $\psi$ realizing $\varphi$.

Next suppose that $2|m_1m_2$. In fact in this case the class $x$ defined above
is fixed $\mod 2$ by every isometry $\varphi$ as above which satisfies $\varphi
(L) = L$: Since $m_1m_2|\beta$ and $c\cong \beta ^2/2\mod 2$, it follows that
$\varphi (x) \equiv x \mod 2$.
Now let $x'$ be a class with $x'\cdot \kappa \equiv 1\mod 2$ and
$(x')^2 \equiv x^2 \mod 4$. We may assume that $x'\neq x$. First consider the
case where $x' = x + b\kappa + \alpha$ and $b\equiv 0 \mod 2$. Thus we may
replace $x'$ by $x+\alpha$. By assumption $\alpha ^2 \equiv 0 \mod 4$. We may
assume that $\alpha $ is primitive (otherwise $\alpha \equiv 0\mod 2$ or
$\alpha$ is congruent to a primitive nonzero element mod 2).  Replacing $\alpha
$ by $\alpha + 2\beta$, where $\beta \in W^\perp$,
replaces $\alpha ^2$ by $\alpha ^2 + 4(\alpha \cdot \beta) + 4\beta ^2$. Since
$\alpha$ is primitive, it is easy to see that there is a choice of $\beta$ so
that $(\alpha + 2\beta)^2 = 0$. Thus we may assume that $\alpha$ is
primitive and that $\alpha ^2=0$. The group  $SO(W)$ includes into the
automorphism group of $L$, and every element of $SO(W)^*$, the set all elements
of $SO(W)$ with spinor norm one, is realized by a diffeomorphism. Moreover an
easy exercise shows that $SO(W)^*$ acts transitively on the set of primitive
$\alpha \in W$ with $\alpha ^2=0$. Thus the set of all possible $x+\alpha$,
with
$\alpha \neq 0$, is contained in a single orbit under the diffeomorphism group.

In case $x' = x + \kappa + \alpha$ with $\alpha ^2
\equiv 2\mod 4$, an argument similar to that
given above shows that we may assume that $\alpha ^2 = 2$ and that every two
classes $x_1 = x + \kappa + \alpha_1$ and $x_2 =  x + \kappa + \alpha_2$ with
$\alpha _i^2 = 2$ are conjugate under the group of diffeomorphisms of $S$ which
fix $\kappa$. Thus there are at most three orbits in this case. \endproof

The following result is really only needed in the case where $m_1m_2$ is even,
since in case $m_1m_2$ is odd we can appeal to (i) of (2.1) above.

\proposition{2.2} Let $S$ be a simply connected elliptic surface and $w$ be a
class in $H^2(S;\Zee/2\Zee)$ with $w\cdot \kappa = 1$. Then after replacing
$S$ with a deformation equivalent elliptic surface, we may assume that there is
a  divisor $\Delta$ on $S$ with $\Delta \cdot \kappa =
1$, and such that all singular fibers of $S$  are
irreducible rational curves with a singular ordinary  double point, i.e\. $S$
is
nodal.
\endstatement
\proof Fix a nodal simply connected elliptic surface with a
section $B$ such that $p_g(B) = p_g(S)$. Using [9], $S$ is deformation
equivalent through elliptic surfaces to a logarithmic transform of $B$ at two
smooth fibers, where the multiplicities of the logarithmic transforms are $m_1$
and $m_2$. Fix one such logarithmic transform $S_0$, and let $\psi\: S \to S_0$
be a diffeomorphism preserving the class of the fiber. Using this
diffeomorphism, we shall identify $S$ and $S_0$. Let $\Delta $ be an element in
$H^2(S; \Zee)$ whose mod 2 reduction is $w$ and such that $\Delta \cdot \kappa
=
1$. We shall show that, by further modifying the complex structure on $S$, we
may assume that $\Delta$ is of type $(1,1)$.

Given $\Delta$, we have the image $i_* ([\Delta]) \in H^2(S;\scrO _S)$, where
$i_*$ is the map induced on sheaf cohomology by the inclusion $\Zee \subset
\scrO _S$. The set of all complex structures of an elliptic surface on $S$ for
which the associated Jacobian surface is $B$ and which are locally isomorphic
to
$S$ is a principal homogeneous space over $H^1(\Pee ^1; \Cal B)$, where $\Cal
B$
is the sheaf of local holomorphic cross sections of $B$ ([9] Chapter 1 Theorem
6.7). Moreover there is a surjective map $H^2(S; \scrO _S) \cong H^2(B; \scrO
_B)\to  H^1(\Pee ^1; \Cal B)$ ([9], Chapter 1, Lemma 5.11). Thus given a
cohomology class $\eta \in H^2(S; \scrO _S)$, we can form the associated
surface
$S^\eta$ and consider the element $i_* ^\eta([\Delta]) \in H^2(S^\eta;\scrO
_{S^\eta}) \cong H^2(S;\scrO _S)$.

\lemma{2.3}
In the above notation,
$i_* ^\eta([\Delta]) = i_* ([\Delta]) + m_1m_2\eta$.
\endstatement
\proof This presumably could be proved by a rather involved direct calculation.
For another argument, note that the map $H^2(S; \scrO_S) \to H^2(S; \scrO_S)$
defined by
$$\eta \mapsto i_* ^\eta([\Delta]) - i_* ([\Delta]) - m_1m_2\eta$$
is holomorphic, since it arises from a variation of Hodge structure. The
argument of Lemma 6.13 in Chapter 1 of [9] shows that $i_* ^\eta([\Delta]) -
i_* ([\Delta]) - m_1m_2\eta$ lies in the countable (not necessarily
discrete) subgroup $H^1(\Pee ^1; R^1\pi _*\Zee)$ of $H^2(S; \scrO_S)= H^1(\Pee
^1; R^1\pi _*\scrO _S)$. This is only possible if the image of the map is
contained in a single point, and since the image contains the origin the map is
identically zero.
\endproof

Returning to the proof of (2.2), since $H^2(S;\scrO _S)$ is divisible, there is
a
choice of $\eta$ so that  $i_* ^\eta([\Delta]) = 0$. For the corresponding
complex structure, $[\Delta]$ is then a $(1,1)$ class. \endproof

Finally we shall describe a way to associate new elliptic surfaces to $S$
which generalizes the construction of the Jacobian surface. Suppose that $S$ is
an elliptic surface over $\Pee ^1$. Let $\eta = \operatorname{Spec}k$ be the
generic point of $\Pee ^1$, where $k=k(\Pee ^1)$ is the function field of the
base curve, let $\bar \eta= \operatorname{Spec}\bar k$, where $\bar k
=\overline{k(\Pee ^1)}$ is the algebraic closure of $k$, and let $S_\eta$ and
$S_{\bar \eta}$ be the restrictions of $S$ to $\eta$ and $\bar \eta$. Thus
$S_\eta$ is a curve of genus one over $k$.

Given an algebraic elliptic surface with a section $\pi \:B\to \Pee ^1$, it has
an associated  Weil-Chatelet group $WC(B)$ [3], which classifies all
algebraic elliptic surfaces $S$ whose Jacobian surface is $B$. As above we let
$B_\eta$ be the elliptic curve over $k$ defined by  the
generic fiber of $B$. By definition $WC(B)$ is the Galois cohomology group
$H^1(G, B_\eta(\bar k))$, where $G = \operatorname{Gal}(\bar k/k)$ and
$B_\eta(\bar k)$ is the group of points of the elliptic curve $B_\eta$ defined
over $\bar k$. There is an exact sequence  $$0 \to \Shah (B) \to WC(B) \to
\bigoplus _{t\in C}H_1(\pi ^{-1}(t); \Bbb Q/\Zee) \to 0.$$  The subgroup $\Shah
(B)$ corresponds to elliptic surfaces without multiple fibers whose Jacobian
surface is isomorphic to $B$, and the quotient describes the possible local
forms for the multiple fibers. Thus if $\xi\in WC(B)$ corresponds to the
surface
$S$, then $S$ has a multiple fiber of multiplicity $m$ at $t\in \Pee ^1$ if and
only if the projection of $\xi$ to $H_1(\pi ^{-1}(t); \Bbb Q/\Zee)$ has order
$m$.

The surface $S$ is specified by an element $\xi$ of $WC(J(S))$, where $J(S)$ is
the Jacobian surface associated to $S$. Let
us recall the recipe for $\xi$ [18]: we have the curve $S_\eta$ and its
Jacobian
$J(S_\eta )$ defined over $k$. The curve $S_\eta$ is a
principal homogeneous space over $J(S_\eta )$, and thus defines a class $\xi
\in
WC(J(S))$, by the following rule: let $\sigma$ be a point of $S_{\bar \eta}$.
Given $g\in \operatorname{Gal}(\bar k/k)$, the divisor $g(\sigma) - \sigma$ has
degree zero on $S_{\bar \eta}$ and so defines an element of $J(S_{\bar \eta})$,
which is easily checked to be a 1-cocycle. The induced cohomology class is
$\xi$.

For every integer $d$ there is an algebraic elliptic surface $J^d(S)$, whose
restriction to the generic fiber $\eta$ is the Picard scheme of divisors of
degree $d$ on the curve $S_\eta$. Thus $J^0(S)
= J(S)$ and $J^1(S) = S$. We claim that, if $S$ corresponds to
the class $\xi \in WC(J(S))$, then $J^d(S)$ corresponds to the class $d\xi$.
Indeed, using the above notation, if $\sigma$ defines a point of $S_{\bar
\eta}$, then $d\sigma$ is a point of $J^d(S_{\bar \eta})$. Thus, the
corresponding cohomology class is represented by $d(g(\sigma) - \sigma)$ and
so is equal to $d\xi$. In particular, if $S$ has a multiple fiber of
multiplicity
$m$ at $t$, then $J^d(S)$ has a multiple fiber of multiplicity $m/\gcd (m, d)$.
Of course if $m|d$ then the multiplicity is one. Finally note that $J(S)$ is
the Jacobian surface of $J^d(S)$ for every $d$ and that $p_g(J^d(S)) = p_g(S)$.

Ideally we would like there to be a Poincar\'e line bundle $\Cal P_d$ over
$S\times _{\Pee ^1}J^d(S)$ such that the restriction of $\Cal P_d$ to the slice
$S\times _{\Pee ^1}\{\lambda\}$ is the line bundle of degree $d$
on the fiber of $S$ over $\pi (\lambda)$ corresponding to $\lambda$. In general
this is too much to ask. However such a bundle exists locally around every
smooth nonmultiple fiber: if $X$ is the inverse image in $S$ of a small disk
$D$ in $\Pee ^1$ such that all fibers on $X$   are smooth and nonmultiple, and
$X_d$ is the corresponding preimage in $J^d(S)$, then there is a Poincar\'e
line
bundle over $X\times _DX_d$. There is also an analogous statement where we
replace a small classical open set in $\Pee ^1$ with an \'etale open set. The
proof for this result is essentially contained in the proof of Theorem 1.3 of
Chapter 7 in [9]. Another construction is given in Section 7 of Part III of
this
series.

\section{3. Suitable line bundles.}

Suppose that we are given a class $w\in H^2(S; \Zee/2\Zee)$ with $w\cdot
\kappa \bmod 2 = 1$ and an integer $p$ with $w^2 \equiv p \mod 4$. Choose once
and for all a complex structure on $S$ for which there is a divisor $\Delta $
with $w= \Delta \bmod 2$. Let $c$ be the integer $(\Delta ^2 - p)/4$. The
principal $SO(3)$-bundle $P$ over $S$ with $w_2(P) = w$ and $p_1(P) = p$ lifts
uniquely to a principal $U(2)$-bundle $P'$ over $S$ with $c_1(P') = \Delta$ and
$c_2(P') = c$. Moreover, by Donaldson's theorem, if $g$ is a Hodge metric on
$S$
corresponding to the ample line bundle $L$, we can identify the moduli space of
gauge equivalence classes of $g$-anti-self-dual connections on $P$ with the
moduli space of $L$-stable rank two vector bundles $V$ over $S$ with $c_1(V) =
\Delta$ and $c_2(V) = c$.

We shall also have to make a choice of the ample line bundle $L$. If $p_g(S) >
0$, then the resulting Donaldson polynomial invariant does not depend on the
choice of $L$, whereas if $p_g(S) = 0$, then the invariant depends on the
chamber containing $c_1(L)$ [11], [12]. We then make the following
definition [17]:

\definition{Definition 3.1} A {\sl wall of type $(\Delta, c)$} is a class
$\zeta \in H^2(S; \Zee)$ such that $\zeta \equiv \Delta \mod 2$ and
$$\Delta ^2 - 4c \leq \zeta ^2 <0.$$
In particular there are no such walls unless  $\Delta ^2 - 4c <0$. Clearly
this definition depends only on $\Delta \bmod 2 = w$ and $p= \Delta ^2 - 4c$,
and
we shall also refer to walls of type $(w, p)$.

Now suppose that $p_g(S) = 0$, i.e\. that $b_2^+(S) = 1$. Let
$$\Omega _S = \{
x\in H^2(S; \Ar): x^2 >0\}.$$
 Let $W^\zeta = \Omega _S \cap (\zeta)^\perp$.  A
{\sl chamber of type $(\Delta, c)$} (or of type $(w, p)$) for $S$  is a
connected
component of the set
$$\Omega _S - \bigcup\{W^\zeta: \zeta {\text{ is a wall of type $(\Delta, c')$,
$c'\leq c$}}\, \}.$$
\enddefinition

For the purposes of algebraic geometry, walls of type $(\Delta, c)$ arise as
follows: let $L$ be an ample line bundle and let $V$ be a rank two bundle over
$S$ with $c_1(V) = \Delta$ and $c_2(V) = c$ which is strictly $L$-semistable.
Let $\scrO _S(F)$ be a destabilizing sub-line bundle. Thus there is an exact
sequence $$0 \to \scrO _S(F) \to V \to \scrO _S(-F+\Delta )\otimes I_Z
\to 0,$$ where $I_Z$ is the ideal sheaf of a codimension two local complete
intersection subscheme. Thus
$$c_2(V) = c = -F^2 + F \cdot \Delta + \ell (Z).$$
Since $\ell(Z)$ is nonnegative, we can rewrite this as
$$ -F^2 + F \cdot \Delta \leq c.$$
Moreover
$$(2F -\Delta)^2 = -4(-F^2 + F \cdot \Delta) + \Delta ^2 \leq \Delta ^2
-4c,$$
so that we can rewrite the last condition by
$$ \Delta ^2 - 4c \leq (2F -\Delta)^2.$$
Using the fact that $L\cdot F = (L\cdot \Delta)/2$, we have
$$L\cdot (2F -\Delta) = 0,$$
and so by the Hodge index theorem $(2F -\Delta)^2 \leq 0$, with equality
holding if and only if $2F -\Delta= 0$ (recall that $S$ is simply connected).
This case cannot arise for us since $\Delta \cdot \kappa = 1$ and thus $\Delta$
is primitive. In particular $\zeta = 2F-\Delta$ is a wall of type $(\Delta,
c)$.
Of course, it is also the cohomology class of a divisor, and thus has type
$(1,1)$. It then follows easily that, if $L_1$ and $L_2$ are two ample line
bundles such that $c_1(L_1)$ and $c_1(L_2)$ lie in the interior of the same
chamber of type $(\Delta, c)$, then a rank two vector bundle $V$ with $c_1(V) =
\Delta$ and $c_2(V)=c$ is $L_1$-stable if and only if it is $L_2$-stable.

With this said, we can make the following definition:

\definition{Definition 3.2}
Let $c$ be an integer, and set $w = \Delta \bmod 2$ and $p = \Delta ^2 - 4c$.
An ample line bundle $L$ is {\sl $(\Delta, c)$-suitable\/} or
{\sl $(w,p)$-suitable\/} if, for all walls $\zeta$ of type $(\Delta, c)$ which
are the classes of divisors on $S$, we have
$\operatorname{sign} f\cdot \zeta = \operatorname{sign} L\cdot \zeta$.
\enddefinition
\medskip

\noindent {\bf Remark.} 1) Suppose that $\zeta$ is a $(1,1)$ class satisfying
$\zeta ^2 \geq 0$. It follows from the Hodge index theorem that if $\zeta \cdot
f>0$, then $\zeta \cdot L >0$ as well. Thus we can drop the requirement that
$\zeta ^2 <0$.

2) In our case $\zeta \equiv \Delta \bmod 2$ and thus $\zeta \cdot \kappa
\equiv 1 \mod 2$. It follows that $\zeta \cdot \kappa \neq 0$ and thus that
$\zeta \cdot f \neq 0$. Thus the condition $\zeta \cdot f \neq 0$ (which was
included as part of the definition in [9]) is always satisfied in our case.

3) In case $b_2^+(S)=1$, $L$ is $(\Delta, c)$-suitable if and only if the class
$\kappa$ lies in the closure of the chamber containing $c_1(L)$.
\medskip

\lemma{3.3} For every $c$, $(\Delta, c)$-suitable ample line bundles exist.
\endstatement
\proof Let $L_0$ an ample line bundle. For $n\geq 0$, let $L_n = L_0 \otimes
\scrO _S(nf)$. It follows from the Nakai-Moishezon criterion that $L_n$ is
ample as well. We claim that if $n> -p (L_0 \cdot f)/2$, then $L_n$ is
$(\Delta, c)$-suitable.

To see this, let $\zeta = 2F - \Delta$ be a wall of type $(\Delta, c)$ with
$$\Delta ^2 - 4c \leq \zeta^2 <0.$$
We may assume that $a = \zeta \cdot f >0$, and must show that $\zeta \cdot L_n
>
0$ as well. The class $ac_1(L_0) - (L_0\cdot f) \zeta$ is perpendicular to $f$.
Since $f^2=0$, we may apply the Hodge index theorem to conclude:
$$0 \geq ac_1(L_0) - (L_0\cdot f) \zeta = a^2L_0^2 -2a(L_0\cdot f) (L_0 \cdot
\zeta) + (L_0\cdot f)^2\zeta ^2.$$
Using the fact that $\zeta ^2 \geq\Delta ^2 - 4c = p $, we find that
$$L_0 \cdot \zeta \geq \frac{a(L_0^2)}{2(L_0\cdot f)} + \frac{\zeta ^2}{2a}(L_0
\cdot f) > \frac{p}{2a}(L_0 \cdot f).$$
Thus
$$\align
L_n \cdot \zeta &= (L_0 \cdot \zeta) + n(\zeta \cdot f) >\frac{p}{2a}(L_0 \cdot
f) -\frac{pa}{2}(L_0\cdot f) \\
&= -\frac{p }{2}(L_0\cdot f)\Bigl(a-\frac{1}{a}\Bigr) \geq 0.
\endalign$$
Thus $L_n$ is $(\Delta, c)$-suitable.
\endproof

In the case where $b_2^+ (S) =1$, we have the following interpretation of
$(\Delta, c)$-suitability.

\lemma{3.4} Suppose that $p_g(S) = 0$. If $L_1$ and $L_2$ are both $(\Delta,
c)$-suitable, then $c_1(L_1)$ and $c_1(L_2)$ lie in the same chamber of type
$(\Delta, c)$. Thus there is a unique chamber $\Cal C(w,p)$ of type
$(\Delta, c)$ which contains the first Chern classes of $(\Delta, c)$-suitable
ample line bundles. Conversely, if $L$ is ample and $c_1(L) \in \Cal C(w,p)$,
then $L$ is $(\Delta, c)$-suitable.
\endstatement
\proof Let $L_1$ and $L_2$ be $(\Delta,
c)$-suitable. Since $p_g(S) =0$, every cohomology class is of type $(1,1)$.
Thus
if $\zeta$ is a wall of type $(\Delta, c)$, then
$$\operatorname{sign} L_1\cdot \zeta =\operatorname{sign} f\cdot \zeta =
\operatorname{sign} L_2\cdot \zeta.$$
This exactly implies that $c_1(L_1)$ and $c_1(L_2)$ are not separated by any
wall $(\zeta )^\perp$.

Conversely suppose that $c_1(L) \in \Cal C(w,p)$, where $\Cal C(w,p)$ is the
unique chamber containing the first Chern classes of $(\Delta, c)$-suitable
ample line bundles. This means in particular that $L\cdot \zeta \neq 0$ for
every
$\zeta$ of type $(\Delta, c)$. The proof of (3.4) shows that $c_1(L) + N[f] \in
\Cal C(w,p)$ for all sufficiently large $N$. Thus, for all $\zeta$ of type
$(\Delta, c)$,
$$\operatorname{sign} L\cdot \zeta = \operatorname{sign} L\cdot \zeta +
N\operatorname{sign} f\cdot \zeta.$$
Since $f\cdot \zeta \neq 0$,
$\operatorname{sign} L\cdot \zeta + N\operatorname{sign} f\cdot \zeta =
\operatorname{sign} f\cdot \zeta$ for all $N\gg 0$. Thus $\operatorname{sign}
f\cdot \zeta = \operatorname{sign} L\cdot \zeta$, and $L$ is $(\Delta,
c)$-suitable. \endproof

\lemma{3.5} Suppose that $p_g(S) = 0$. The chamber $\Cal C(w,p)$ is the unique
chamber of type $(w,p)$ which contains $\kappa$ in its closure. Thus every
diffeomorphism $\psi$ of $S$ satisfies $\psi ^*\Cal C(w,p) = \pm \Cal C(\psi
^*w,p)$. More generally, if $S$ and $S'$ are two elliptic surfaces with $p_g=0$
and $\psi \: S \to S'$ is a diffeomorphism, then $\psi ^*\Cal C(w,p) =  \pm
\Cal
C(\psi ^*w,p)$.
\endstatement
\proof Let $\Cal C_1$ and $\Cal C_2$ be two distinct chambers which contain
$\kappa$ in their closures. Let $\zeta$ be a wall separating $\Cal C_1$ and
$\Cal C_2$. We may assume that $\zeta \cdot x > 0 $ for all $x\in \Cal C_1$ and
$\zeta \cdot x  < 0$ for all $x \in \Cal C_2$. Thus $0\leq\zeta \cdot \kappa
\leq
0$, so that $\zeta \cdot \kappa = 0$. However this contradicts the fact that
$\zeta \cdot \kappa \neq 0$. Thus there is at most one
chamber containing $\kappa$ in its closure. We have seen in the proof of Lemma
3.3 that, for all ample line bundles $L$ and integers $N\gg 0$, $c_1(L) +
N\kappa \in \Cal C(w,p)$. Thus $\kappa + (1/N)c_1(L) \in \Cal C(w,p)$. It
follows that $\kappa$ indeed lies in the closure of $\Cal C(w,p)$, so that
$\Cal C(w,p)$ is the unique chamber with this property.

To see the final statement, we use [8] to see that every diffeomorphism
$\psi$ of $S$ satusfies $\psi ^*\kappa = \pm \kappa$. Thus $\pm \kappa$ lies in
the closure of $\psi ^*\Cal C(w,p)$. Clearly, if $\Cal C$ is a chamber of type
$(w,p)$, then $\psi ^*\Cal C$ is a chamber of type $(\psi ^*w,p)$. It follows
that $\psi ^*\Cal C(w,p) =  \pm \Cal
C(\psi ^*w,p)$. The statement about two different surfaces is proved similarly.
\endproof

\definition{Definition 3.6} The chamber described in Lemma 3.4 will be called
the {\sl suitable chamber\/} of type $(\Delta, c)$ or of type $(w,p)$ or the
{\sl $(\Delta, c)$-suitable\/} or {\sl $(w, p)$-suitable chamber\/}.
\enddefinition

\section{4. The geometric meaning of suitability.}

The goal of this section is to describe the meaning of $(\Delta,
c)$-suitability. Given the bundle $V$ on $S$, it defines by restriction a
bundle $V|f$ on each fiber $f$. Our main result says essentially that $V$ is
stable for one, or equivalently all,  $(\Delta,
c)$-suitable line bundles $L$ if and only if $V|f$ is semistable for almost all
$f$.

It will be more convenient to use the language of schemes to state this result.
As in Section 2, let $k(\Pee ^1)$ denote the function field of $\Pee ^1$ and
let
$\overline  {k(\Pee ^1)}$ be the algebraic closure of $k(\Pee ^1)$. Set $\eta =
\operatorname{Spec}k(\Pee ^1)$ and $\bar\eta = \operatorname{Spec}
\overline{k(\Pee ^1)}$. Thus $\eta$ is the generic point of $\Pee ^1$. Let
$S_\eta= S\times _{\Pee ^1}\eta$ be the generic fiber of $\pi$ and let $S_{\bar
\eta} = S\times _{\Pee ^1}\bar\eta$. Here $S_\eta$ is a curve of genus one over
the field $k(\Pee ^1)$ and $S_{\bar\eta}$ is the curve over $\overline{k(\Pee
^1)}$ defined by extending scalars. Let $V_\eta$ and $V_{\bar \eta}$ be the
vector bundles over $S_\eta$ and $S_{\bar\eta}$ respectively obtained by
restricting $V$. We can then define stability and semistability for $V_{\bar
\eta}$ and $V_\eta$; for $V_\eta$, a destabilizing subbundle must also be
defined over $k(\Pee ^1)$. Trivially, if $V_\eta$ is unstable (resp\. not
stable)
then $V_{\bar \eta}$ is unstable (resp\. not stable). Thus if $V_{\bar \eta}$
is stable, then $V_\eta$ is stable as well.

\lemma{4.1}  $V_\eta$ is semistable if and only if $V_{\bar \eta}$ is
semistable.
\endstatement
\proof We have seen that, if $V_\eta$ is  not semistable, then $V_{\bar \eta}$
is not semistable. Conversely suppose that $V_{\bar \eta}$ is not semistable.
Then there is a canonically defined maximal destabilizing line subbundle of
$V_{\bar \eta}$, which thus is fixed under by every element of
$\operatorname{Gal}(\overline  {k(\Pee ^1)}/k(\Pee ^1)$. By standard descent
theory  this line subbundle must then be defined over $k(\Pee ^1)$. Thus
$V_\eta$ is  not semistable.
\endproof

\noindent {\bf Remark.} If
$V_{\bar \eta}$ is strictly semistable, it is typically the case that $V_\eta$
is actually stable.
\medskip

\lemma{4.2} In case $\Delta \cdot \kappa = 1$, the
bundle $V_\eta$ is semistable if and only if it is stable.
\endstatement
\proof
First assume that $m_1m_1 \equiv 1\mod 2$. In this case $V_\eta$ has odd
fiber degree, and so there are no strictly semistable bundles over
$\overline{k(\Pee ^1)}$. Hence, if $V_\eta$ is semistable, then by (4.1)
$V_{\bar \eta}$ is semistable and therefore stable. Thus $V_\eta$ is
stable by the remarks preceding (4.1).

In case $m_1m_1 \equiv 0\mod 2$, suppose that  $V_\eta$ is strictly semistable.
Thus there  is a line bundle on $S_\eta$ of degree $m_1m_2/2$. There would thus
exist a divisor $D$ on $S$ with $D\cdot f = m_1m_2/2$. Since $f =
m_1m_2\kappa$, this possibility cannot occur. Thus $V_\eta$ is stable.
\endproof

Here then is the theorem of this section:

\theorem{4.3} Let $V$ be a rank two vector bundle on $S$ with $c_1(V) =
\Delta$ and $c_2(V) = c$ and let $L$ be a $(\Delta, c)$-suitable ample line
bundle. Then $V$ is $L$-stable if and only if the restriction $V_\eta$ of $V$
to the generic fiber $S_\eta$ is stable.
\endstatement
\proof First suppose that $V$ is $L$-stable. Let $F_\eta$ be a subbundle of
$V_\eta$ of rank one. Then there is a divisor $F$ on $S$ such that $\scrO
_S(F)$ restricts to $F_\eta$ and an inclusion $\scrO _S(F) \to V$. Hence
there is an effective divisor $D$ and an inclusion
$\scrO_S(F+D) \to V$ and the cokernel is torsion free. Since $F_\eta$ is a
subbundle of $V_\eta$, the divisor $D$ cannot have positive intersection number
with $f$. As $D$ is effective it is supported in the fibers of $\pi$ and so
$F$ and $F+D$ have the same restriction to the generic fiber. We may thus
replace $F$ by $F+D$.  Then $V/\scrO _S(F)$ is torsion free. Hence there is an
exact sequence  $$0 \to \scrO _S(F) \to V \to \scrO _S(\Delta - F)\otimes I_Z
\to
0,$$
where $Z$ is a codimension two subscheme of $S$. Thus  $$\Delta ^2 - 4c \leq
(2F-\Delta )^2.$$ Since $V$ is $L$-stable, $L\cdot (2F-\Delta ) <0$. It follows
from Definition 3.2 and 2) of the remark following it that $f\cdot (2F-\Delta )
<0$ as well. Thus $\deg F_\eta <\deg V_\eta /2$, which says that $V_\eta$ is
stable.

Conversely suppose that $V_\eta$ is stable. Let $\scrO _S(F)$ be a sub-line
bundle of $V$, where we may assume that $V/\scrO _S(F)$ is torsion free.
Reversing the argument above shows that $f\cdot (2F-\Delta ) <0$ and therefore
that $L\cdot (2F-\Delta )<0$ as well. Thus $V$ is $L$-stable.
\endproof

\corollary{4.4}  Let $V$ be a rank two vector bundle on $S$ with $c_1(V) =
\Delta$ and $c_2(V) = c$. Then the  following are equivalent:
\roster
\item"{(i)}" There exists a $(\Delta, c)$-suitable ample line
bundle $L$ such that $V$ is $L$-stable.
\item"{(ii)}" $V$ is $L$-stable for every $(\Delta, c)$-suitable ample line
bundle $L$.
\item"{(iii)}" $V_\eta$ is stable.
\item"{(iv)}" $V_{\bar \eta}$ is semistable.
\item"{(v)}" The restriction $V|\pi^{-1}(t)$ is semistable for almost all $t
\in
\Pee ^1$.
\item"{(vi)}" There exists a $t\in \Pee^1$ such that $\pi^{-1}(t)$ is smooth
and the  restriction $V|\pi^{-1}(t)$ is semistable.
\endroster
\endstatement
\proof By (4.1) and (4.2), (iii) and (iv) are equivalent, and by (4.3) (i)
$\implies$  (iii) $\implies$ (ii). The implication (ii) $\implies$ (i) is
trivial. The implication (iv) $\implies$ (v) follows from the openness of
semistability in the Zariski topology in the sense of schemes, and the
implication (v) $\implies$ (vi) is trivial. To see that (vi) $\implies$ (iv),
suppose that $V_{\bar \eta}$ is not semistable. Then a destabilizing sub-line
bundle extends to give a sub-line  bundle over the pullback of $S$ to some
finite base change of $\Pee^1$. Thus $V|\pi^{-1}(t)$ is unstable for every
$t\in
\Pee^1$ such that $\pi^{-1}(t)$ is smooth, and so (vi) $\implies$ (iv).
\endproof

\noindent {\bf Remark.} In case $\Delta \cdot \kappa \equiv 0\mod 2$, there can
exist strictly semistable bundles on $S_\eta$ of degree $\Delta \cdot f$. There
are examples of rank two bundles $V$ on $S$ with $c_1(V) = \Delta$ and $V_\eta$
strictly semistable such that $V$ is either stable, strictly semistable, or
unstable (cf\. [8]).

\section{5. Donaldson polynomials and the main theorems.}

As above we let $S$ denote a simply connected elliptic surface with $p_g(S)
\geq 0$. Fix $w = \Delta \bmod 2$ and let $p$ be an integer satisfying $w^2
\equiv p \bmod 4$. For $p_g(S) >0$, there is the Donaldson polynomial $\gamma
_{w,p}(S)$ corresponding to the $SO(3)$-bundle $P$ with invariants $w$ and $p$.
Here for simplicity we shall always choose the orientation
on the moduli space which agrees with the natural complex orientation.
The polynomial $\gamma_{w,p}(S)$ is invariant up to sign under
self-diffeomorphisms $\psi$ of  $S$  such that $\psi ^*w=w$.  If  $p_g(S) = 0$,
then we have the distinguished chamber $\Cal C(w,p)$ which contains $\kappa$ in
its closure. We shall then use $\gamma_{w,p}(S)$ to denote  the Donaldson
polynomial for $S$ with respect to the chamber $\Cal C(w,p)$, again with the
orientation chosen to be the complex orientation. Since $\psi ^*\Cal
C(w,p) = \pm \Cal C(\psi ^*w,p)$, the invariant $\gamma_{w,p}(S)$ is again
natural up to sign under orientation-preserving self-diffeomorphisms which fix
$w$. Of course, there are only finitely many choices for $w$, so that there is
a subgroup of finite index in the full group of diffeomorphisms fixing $[f]$
which will also fix $w$.

\lemma{5.1} For every choice of $w$ and $p$, $\gamma_{w,p}(S)$ lies in $\Bbb
Q[q_S, \kappa _S]$. Moreover, if for some choice of $w$ and $p$,
$\gamma_{w,p}(S)$ does not lie in $\Bbb Q[q_S]$, then every diffeomorphism
$\psi$ from $S$ to another simply connected elliptic surface $S'$ satisfies
$\psi
^*\kappa _{S'} = \pm \kappa _S$.
\endstatement
\proof The set of automorphisms of $H_2(S; \Zee)$ of the form $\psi _*$, where
$\psi$ is a diffeomorphism satisfying $\psi _*([f]) = [f]$, $\psi ^*w = w$, and
$\psi ^* \gamma_{w,p}(S) = \gamma_{w,p}(S)$, is a subgroup of finite index in
the group of all isometries of $H_2(S;\Zee)$ preserving $[f]$, by [8] Part I
Theorem 6 and [9] Chapter 2 Theorem 6.5. Thus by [9] Chapter 6 Theorem 2.12,
$\gamma_{w,p}(S) \in \Bbb Q[q_S, \kappa _S]$. Moreover $\kappa$ is the unique
such class. The last statement of the lemma is then clear.
\endproof

Next let us discuss the effect of blowing up. Suppose that $\rho\: \tilde S\to
S$ is the $r$-fold blowup of $S$, and let the exceptional classes in
$H_2(\tilde
S)$ be denoted by $e_1, \dots, e_r$. Likewise let $S'$ be another simply
connected elliptic surface and let $\rho'\:\tilde S'\to S'$ be the $r$-fold
blowup of $S'$, with exceptional classes $e_1', \dots, e_r'$. If $\psi \:
\tilde
S \to \tilde S'$ is a diffeomorphism, then $\psi ^*e_i' = \pm e_j$ for a
uniquely determined $j$, by [8] Part I Theorem 7 and [9] Chapter 6 Corollary
3.8.
It follows that, if $w'\in H^2(\tilde S'; \Zee/2\Zee)$ is of the form $(\rho
')^*w_0'$ for some $w_0' \in H^2(S'; \Zee/2\Zee)$, then there is a $w_0 \in
H^2(S; \Zee/2\Zee)$,
such that $\psi ^*w' = \rho ^*w_0$.  Finally, we shall need the following
extension of [9] Chapter 6 Theorem 3.1:

\proposition{5.2} Let $w_0 \in H^2(S; \Zee/2\Zee)$ and let $\rho \:\tilde S\to
S$ be the $r$-fold blowup of $S$. If $b_2^+(S) =1$, assume moreover that
$\gamma _{\rho ^*w_0, p}$ is defined with respect to some chamber $\Cal D$. Let
$\Cal C$ be a chamber of type $(w_0, p)$ on $H^2(S; \Bbb R)$ such that $\Cal D$
contains $\rho ^*\Cal C$ in its closure. Then  $$\gamma _{\rho ^*w_0,
p}|\rho ^*H_2(S; \Zee) = \gamma _{w_0, p},$$  where if $b_2^+(S)=1$, the
polynomial $\gamma _{w_0, p}$ is defined with respect to the chamber $\Cal C$.
\endstatement
\medskip

Here, in case $p_g(S) = 0$, the chamber $\Cal C$ does not in general determine
a unique chamber $\Cal D$ on $\tilde S$. However the conclusion of the
proposition implies in particular that the value of $\gamma _{\rho ^*w_0, p}$
on
classes in $\rho ^*H_2(S; \Zee)$ is independent of the chamber for $\tilde S$
of
type $(\rho ^*w_0, p)$ which contains $\Cal C$ in its closure.

This result follows from standard gauge theory techniques [5]. It can also
be proved in our case via algebraic geometry, using the blowup formulas for
instance in [8]. Since it does not appear with an explicit proof in the
literature, we shall outline a proof in the only case that concerns us, where
the chamber $\Cal C$ contains the first Chern class of an
ample line bundle. We shall just write down the argument in the most
interesting case, where $p_g(S) = 0$. We shall also assume that the moduli
space  has the  expected dimension, although the arguments given here can
easily
be extended to handle the general case.

By induction we may assume that $\rho \: \tilde S \to S$ is the blowup of $S$
at
a single point $p$. Let $E$ be the exceptional curve and $e$ be its cohomology
class. We shall usually identify $H^2(S)$ with its image in $H^2(\tilde S)$
under
$\rho ^*$. Let $\Cal D$ be a chamber for $\tilde S$ of type $(\rho ^*w_0, p)$
containing $\Cal C$ in its closure and let $\zeta$ be a wall for $\Cal D$. Then
$\zeta = \zeta ' + ae$, where $\zeta ' \in H^2(S; \Zee)$ and $a\in \Zee$ (in
fact $2|a$ since $\zeta \equiv \Delta \mod 2$). After possibly reflecting in
$e$, which is realized by an orientation-preserving diffeomorphism $r_e$ of
$\tilde S$, we may assume that $a\geq 0$: Indeed, $r_e^*$ switches the two
possible chambers corresponding to $\pm a$, and so, if $\gamma _1$ and $\gamma
_2$ are the two invariants corresponding to the two choices of chambers, then
$r_e^*\gamma _1 = \gamma _2$. Since  $r_e^*|H_2(S; \Zee)$ is the identity, it
suffices to prove the result for either chamber. So we can assume that $a\geq
0$.

Since $\Cal C$ is
in the closure of $\Cal D$, if $x\in \Cal C$, then $x\cdot \zeta '= x\cdot
\zeta
\geq 0$. Conversely, if we start with an ample line bundle $L$ on $S$ such that
$c_1(L) \in \Cal C$, then for all $N\gg 0$, $Nc_1(L) - e$ is the first Chern
class of an ample line bundle $L_N$ on $\tilde S$. Moreover $(Nc_1(L) - e)\cdot
\zeta \geq N(c_1(L) \cdot \zeta ') \geq 0$. It follows from this that
$c_1(L_N)$
lies in $\Cal D$, and if $c_1(L)$ is in the interior of $\Cal C$ then
$c_1(L_N)$
is in the interior of $\Cal D$.

Consider rank two vector bundles $\tilde V$ on $\tilde S$ with $c_1(\tilde V) =
\rho ^*\Delta$ and $c_2(\tilde V) = c$. Set $V = (\rho
_*\tilde V)\spcheck{}\spcheck$. Then $V$ is a rank two vector bundle on $S$
with
$c_1(V) = \Delta$ and $c_2(V) \leq c_2(\tilde V)$, where equality holds if and
only if $\tilde V = \rho ^*V$. The arguments of the proof of
Theorem 5.5 in Part II of [8], which essentially just depend on the
determinant of $\tilde V$ being a pullback, show the following. There is a
constant $N_0$, depending only on $L$ and $c$, such that, for all $N\geq N_0$,
if $\tilde V$ is $L_N$-stable then $V$ is $L$-semistable, and conversely if $V$
is $L$-stable then $\tilde V$ is $L_N$-stable. Moreover the map $V\mapsto \rho
^*V$ defines an open immersion of schemes from the moduli space of $L$-stable
rank two vector bundles on $S$ with $c_1=\Delta$ and $c_2=c$ to the
corresponding moduli space for $\tilde S$ and $c_1 = \rho ^*\Delta$.

To evaluate the Donaldson polynomial on $\tilde S$ on a collection of classes
of the form $\rho ^*\alpha$, represent $\alpha$ by a smoothly embedded Riemann
surface $C$ on $S$ which does not pass through $p$, the center of the blowup,
and choose a theta characteristic on $C$. This choice leads to a divisor $D_C$
on the moduli space. By definition $\tilde V$ lies in $D _C$ if and only if the
Dirac operator coupled to the ASD connection induced on $C$ has a kernel. From
this it is clear that $\tilde V$ lies in $D _C$ if and only if $V$ lies in the
corresponding divisor on the moduli space for $S$ of bundles with $c_1 =
\Delta$
and $c_2 = c_2(V)\leq c$. An easy counting argument then shows that, if $d$ is
the dimension of the moduli space for $\tilde S$ and we choose $C_1, \dots,
C_d$
in general position and general theta characteristics on $C_i$, then $\tilde V$
lies in the intersection $D _{C_1} \cap \dots \cap D_{C_d}$ if and only if
$\tilde V = \rho ^*V$ and $V$ lies in the corresponding intersection for the
moduli space of $S$. As $\#(D _{C_1} \cap \dots \cap D_{C_d})$ calculates the
value $\gamma _{\rho ^*w_0, p}([C_1], \dots, [C_d])$, it is then clear that
$$\gamma _{\rho ^*w_0, p}|\rho ^*H_2(S; \Zee) = \gamma _{w_0, p}. \qed $$
\medskip

Assuming Theorems 1.4 and 1.5, we can now state the main results of this
series of papers.

\theorem{5.3} Suppose that $S$ and $S'$ are two simply connected elliptic
surfaces with $p_g(S) = p_g(S') =1$. Suppose that neither $S$ nor $S'$ is a
$K3$
surface. Let $\tilde S$ and $\tilde S'$ be two blowups of $S$ and $S'$, and let
$\psi \: \tilde S \to \tilde S'$ be a diffeomorphism. Identify $H^2(S; \Zee)$
with its image in $H^2(\tilde S; \Zee)$ under the natural map, and similarly
for
$H^2(S'; \Zee)$. Then $\psi ^*\kappa _{S'} = \pm \kappa _S$.
\endstatement
\proof Arguing as in Corollary 3.6 of Chapter 6 of [9], we see that it
suffices to show that some $\gamma _{w,p}(S)$ actually involves $\kappa _S$. If
$m_1m_2\equiv 1\mod 2$, then the coefficient of $\kappa _S^2$ in $\gamma
_{w,p}(S)$, for the choice of $p$ given in Theorem 1.5 (i) corresponding to the
two-dimensional moduli space, is $2m_1^2m_2^2 -m_1^2 -m_2^2$. This number is
zero
if  $m_1 = m_2 =1$, in which case $S$ is a $K3$ surface. Otherwise $2m_1^2m_2^2
-m_1^2 -m_2^2 >0$. Thus the coefficient  of $\kappa _S^2$ is nonzero.

If $m_1m_2\equiv 0\mod 2$, then the expected dimension of the moduli space is
$4c -\Delta ^2 - 6 \equiv \Delta ^2 \mod 2$. Moreover $\Delta \cdot \kappa _S =
1$ and $K_S = (2m_1m_2 - m_1 - m_2)\kappa$. Since exactly one of $m_1, m_2$ is
even, $\Delta \cdot K_S \equiv 1 \mod 2$. Thus by the Wu formula $\Delta ^2
\equiv 1\mod 2$, and the dimension of the moduli space is odd. It follows that
every nonzero invariant must involve $\kappa_S$. Since nonzero invariants exist
by Theorem 1.4 or more generally by Donaldson's theorem on the nonvanishing of
the invariants, there are choices of $p$ for which
$\gamma _{w,p}(S)$ actually involves $\kappa _S$.
\endproof

\theorem{5.4} Suppose that $S$ and $S'$ are two elliptic surfaces with
$p_g(S) = p_g(S') \geq 1$ with finite cyclic fundamental group and multiple
fibers of multiplicities $\{m_1,m_2\}$ and $\{m_1', m_2'\}$, respectively, that
$\tilde S$ and $\tilde S'$ are two blowups of $S$ and $S'$, and that $\tilde S$
and $\tilde S'$ are diffeomorphic. Then $\{m_1,m_2\}=\{m_1', m_2'\}$. Hence $S$
and $S'$ are deformation equivalent.
\endstatement
\proof As in [9] we can reduce to the simply connected case. Using Theorem
1.3, we know that $m_1m_2 = m_1'm_2'$ and that, if  $\psi \: \tilde
S \to \tilde S'$ is a diffeomorphism, then $\psi _*(H_2(S; \Zee)) = H_2(S';
\Zee)$ under the identification of  $H_2(S; \Zee)$ with a subspace of
$H_2(\tilde S;\Zee)$ and likewise for $S'$. Fix a class $w\in H^2(S'; \Zee)$
with $w\cdot \kappa_{S'}=1$. Thus  $$\gamma _{\psi ^*w,p}(\tilde S)|H_2(S;
\Zee)
= \gamma _{w,p}(\tilde S')|H_2(S'; \Zee)$$ under the natural identifications.
Moreover $\psi ^*w\cdot \kappa _S = w\cdot \kappa _{S'} =1$. Thus the
Donaldson polynomial invariants for the minimal surfaces $S$ and $S'$ for the
values $\psi ^*w$ and $w$ respectively are equal. If $m_1m_2 = m_1'm_2'\equiv
1\mod 2$, then  $$(p_g(S) + 1)m_1^2m_2^2 -m_1^2 -m_2^2 = (p_g(S) +
1)(m_1')^2(m_2')^2 -(m_1')^2 -(m_2')^2.$$ Thus $m_1m_2 = m_1'm_2'$ and $m_1^2
+m_2^2= (m_1')^2 +(m_2')^2$. It follows that $(m_1 + m_2)^2 = (m_1' + m_2')^2$
and so $m_1 + m_2 = m_1' + m_2'$. We can thus determine the elementary
symmetric
functions of $m_1$ and $m_2$ from the diffeomorphism type, and hence the
unordered pair.

If $m_1m_2 = m_1'm_2'\equiv 0\mod
2$, then, assuming that $2|m_1$, it follows from Theorem 1.4 that we can
determine $m_1m_2$ and $m_2$. Thus, we can determine $m_1$ as well.
\endproof

\theorem{5.5} Suppose that $S$ and $S'$ are two nonrational
elliptic surfaces with finite cyclic fundamental group and with $p_g(S) =
p_g(S')
=0$ with multiple fibers of multiplicities $\{m_1,m_2\}$ and $\{m_1', m_2'\}$,
respectively, that $\tilde S$ and $\tilde S'$ are two blowups of $S$ and $S'$,
and that $\tilde S$ and $\tilde S'$ are diffeomorphic. Suppose further that
$m_1m_2 \equiv 0 \mod 2$. Then $m_1'm_2' \equiv 0\mod 2$, and
$\{m_1,m_2\}=\{m_1', m_2'\}$.
\endstatement
\proof We may again reduce to the simply connected case. Note that every
diffeomorphism $\psi$ from $\tilde S$ to $\tilde S'$ sends the subspace
$H^2(S';\Zee)$ to $H^2(S; \Zee)$. Choose a class $w\in H^2(S'; \Zee/2\Zee)$
with $w\cdot \kappa _{S'} = 1$. We must have  $\psi ^*\Cal C(w,p)= \pm
\Cal C(\psi ^*w,p)$, by Lemma 3.5. As in the preceding argument, we are
immediately reduced to comparing the Donaldson invariants for the surfaces $S$
and $S'$. Let us first show that  $m_1'm_2' \equiv 0\mod 2$. First note that
the
two-dimensional invariant corresponds to $-p=5> 4=2(4p_g(S) +2)$. Thus we are
in
the stable range and can apply Theorem 1.4 to conclude that the leading
coefficient of $\gamma _{w, -5}(S) = m_2$.  Since $S$ is not rational $m_2>1$.
But if $m_1'm_2'\equiv 1\mod 2$, then by (i) of Theorem 1.5 the leading
coefficient of $\gamma _{w, -5}(S')$ is 1, a contradiction. Hence
$m_1'm_2'\equiv 0\mod 2$ and the leading coefficient of $\gamma _{w,-5}(S')$ is
just $m_2'$. It follows that $m_2 = m_2'$ and, by Bauer's result
(Theorem 1.2), that $(m_1^2-1)(m_2^2 -1) = ((m_1')^2-1)((m_2')^2 -1)$. Thus
$m_1
= m_1'$ as well. \endproof

\theorem{5.6} Suppose that $S$ and $S'$ are two nonrational
elliptic surfaces with finite cyclic fundamental group and with $p_g(S) =
p_g(S')
=0$ with multiple fibers of multiplicities $\{m_1,m_2\}$ and $\{m_1', m_2'\}$,
respectively, that $\tilde S$ and $\tilde S'$ are two blowups of $S$ and $S'$,
and that $\tilde S$ and $\tilde S'$ are diffeomorphic. Suppose further that
$m_1m_2 \equiv 1 \mod 2$. Then $m_1'm_2' \equiv 1\mod 2$, and
$\{m_1,m_2\}=\{m_1', m_2'\}$.
\endstatement
\proof As before we pass to the simply connected case. If $m_1'm_2' \equiv
0\mod 2$, then by (5.5) $m_1m_2 \equiv 0\mod 2$ as well, a contradiction. Thus
$m_1'm_2' \equiv 1\mod 2$. Using (i) and (ii) of Theorem 1.5, we see that the
Donaldson polynomials determine the quantities
$$\align
A &= m_1^2m_2^2 - m_1^2 - m_2^2+1=(m_1^2-1)(m_2^2-1);\\
B &= m_1^4m_2^4 - m_1^4 - m_2^4+1=(m_1^4-1)(m_2^4-1).
\endalign$$
We must show that $A$ and $B$ determine $\{m_1, m_2\}$ provided that both $m_1$
and $m_2$ are greater than one. This is just a matter of elementary algebra:
let $\sigma _1= m_1^2 + m_2^2$ and $\sigma _2 = m_1^2m_2^2$. Then $\sigma _1$
and $\sigma _2$ are the elementary symmetric functions in $m_1^2$ and $m_2^2$
and thus determine $\{m_1^2, m_2^2\}$. As $m_1$ and $m_2$ are positive the
knowledge of $\{m_1^2, m_2^2\}$ determines $\{m_1, m_2\}$.

To read off $\sigma _1$ and $\sigma _2$ from $A$ and $B$, note that if $A\neq
0$ then $$\frac{B}{A} = (m_1^2+1)(m_2^2+1) =  m_1^2m_2^2 + m_1^2 + m_2^2+1.$$
Thus $2\sigma _2 = B/A + A-2$ and $2\sigma _1 = B/A -A$.
provided that $A\neq 0$. Now $A= 0$ if $m_1$ or $m_2$ is one, and otherwise
$A\geq 1$. Thus, provided neither of $m_1$ or $m_2$ is one, $A$ and $B$
determine
$\sigma _2$ and $\sigma _1$.
\endproof

\enddocument